\def\N{\mathbb N}
\def\Z{\mathbb Z}
\def\pf{\begin{proof}}
\def\pfk{\end{proof}}

\documentclass{ita}
\usepackage{amsfonts,amscd}
\usepackage{amsmath}
\usepackage{amsthm}
\usepackage{amssymb,epsfig}
\usepackage{a4}
\usepackage{color}
\usepackage{graphics}
\newtheorem{thm}{Theorem}[section]
\newtheorem{definition}[thm]{Definition}
\newtheorem{lemma}[thm]{Lemma}
\newtheorem{corollary}[thm]{Corollary}

\newtheorem{proposition}[thm]{Proposition}

\begin{document}
\title{Palindromic complexity of infinite words associated
with non-simple Parry numbers}

\author{L$\!$'ubom\'\i ra Balkov\'a}
\author{Zuzana Mas\'akov\'a}
\address{Doppler Institute for Mathematical Physics and Applied
Mathematics \& Department of Mathematics, FNSPE, Czech Technical
University, Trojanova 13, 120 00 Praha 2, Czech Republic\\
\email{l.balkova@centrum.cz \&\ masakova@km1.fjfi.cvut.cz }}



\begin{abstract}
We study the palindromic complexity of infinite words $u_\beta$,
the fixed points of the substitution over a binary alphabet,
$\varphi(0)=0^a1$, $\varphi(1)=0^b1$, with $a-1\geq b\geq 1$,
which are canonically associated with quadratic non-simple Parry
numbers $\beta$.
\end{abstract}

\runningtitle{Palindromic complexity of beta-integers}
\runningauthors{L$\!$'. Balkov\'a, Z. Mas\'akov\'a}

\maketitle

\section{Introduction}

Palindromes, words which stay the same when read backwards, is
a~popular linguistic game. The longest palindrome listed in Oxford
Dictionary is ``tattarrattat''. This word arose in the fantasy of
James Joyce and signified to tap on the door in his novel Ulysses.
Not only words, but also palindromic numbers have been attracting
attention for ages. Such numbers appear already in a~sanskrit
manuscript Ganitasaras{\^a}mgraha dated around 850 AD. However,
palindromes occur in serious disciplines as well. For instance,
nucleotides in most of human genomes form palindromic sequences.

In this paper, we focus on palindromes in infinite aperiodic words
which can serve as models for one dimensional quasicrystals.
Quasicrystals are materials with long-range orientational order
demonstrated by sharp bright spots on their diffraction images,
nevertheless, of non-crystallographic symmetry which testifies
their aperiodicity. For their particular features, these materials
have been in the center of interest of physicists, chemists, and
mathematicians since their discovery in 1982~\cite{Shechtman}. To
understand physical properties of these materials, it is important
to describe their combinatorial properties, such as factor
complexity which corresponds to the number of local configurations
of atoms, but also the palindromic structure which describes local
symmetry of the material. Moreover, the palindromic structure of
the infinite words turns out to be important for the description
of the spectral properties of Schr\H odinger operators with
potentials adapted to aperiodic structures~\cite{HKS}.

Here, we concentrate on infinite words $u_{\beta}$ coding numbers
with integer $\beta$-expan\-sion. Such words are defined over a
finite  alphabet if the base $\beta$ of the numeration system is a
Parry number, i.e.,\ such that its R\'enyi expansion of unity
$d_\beta(1)$ is eventually periodic. For the definition of
$d_\beta(1)$ see Section~\ref{sec:preli}. If $d_\beta(1)$ is
finite, $\beta$ is called a simple Parry number. For such $\beta$,
diverse combinatorial aspects of infinite words $u_\beta$ have
been investigated. In~\cite{FrMaPe1}, one describes the factor
complexity of a large subclass of these words.
In~\cite{AmFrMaPe1}, the palindromic complexity is given for the
same subclass. In paper~\cite{return},the authors determine the
number of return words. The balance properties of quadratic simple
Parry numbers are studied in~\cite{turek}.

Much less is known about infinite words $u_{\beta}$ associated
with non-simple Parry numbers. Some preliminary results about
their factor complexity are given in~\cite{FrMaPe2}. Otherwise,
one finds concrete results only about infinite words $u_\beta$
where $\beta$ is a quadratic non-simple Parry number. In this
case, $u_\beta$ is a fixed point of the substitution
$\varphi(0)=0^{a}1, \ \varphi(1)=0^{b}1, \ a>b \geq 1$. The factor
complexity of such infinite words has been determined
in~\cite{FrMaPe1}. In~\cite{balance}, one studies their
arithmetical and balance properties.

The aim of this paper is to determine the palindromic complexity
of infinite words $u_\beta$ associated with non-simple Parry
numbers. In fact, we explain that the only interesting case is
again the case of quadratic non-simple Parry numbers, since the
language of the infinite word $u_\beta$ for other non-simple Parry
numbers $\beta$ contains only finitely many palindromes.

\section{Notations and definitions}\label{sec:preli}

\subsection{Combinatorics on words} \label{words}

An {\em alphabet} $\cal A$ is a~finite set of symbols called {\em
letters}. A concatenation $w=w_1w_2\cdots w_n$ of letters is
called a {\em word} of length $|w|=n$. The set of all finite
words, including the empty word $\varepsilon$, is denoted by $\cal
A^{*}$, and has the structure of a free monoid with respect to the
operation of concatenation. We deal also with right-sided infinite
words, and, if need be, bidirectional infinite words, i.e., with
sequences of letters $u=u_0u_1u_2\cdots$, or $u=\cdots
u_{-1}u_0u_1u_2\cdots$, respectively. A~finite word $w$ is called
a~{\em factor} of a word $u$ (finite or infinite) if there exist
a~finite word $w^{(1)}$ and a~word $w^{(2)}$ (finite or infinite)
such that $u=w^{(1)}ww^{(2)}$. The word $w$ is a~{\em prefix} of
$u$ if $w^{(1)}=\varepsilon$. Analogously, $w$ is a~{\em suffix}
of $u$ if $w^{(2)}=\varepsilon$. A~concatenation of $k$ words $z$
will be denoted by $z^k$, a~concatenation of infinitely many words
$z$ by $z^{\omega}$. An infinite word $u$ is said to be {\em
eventually periodic} if there exist words $w,z$ such that
$u=wz^{\omega}$. Let $u=u_1u_2u_3\cdots, u_i \in {\cal A}$, then
$u_1^{-1}u=u_2u_3\cdots$. The {\em language} of $u$, denoted by
${\cal L}(u)$, is the set of all factors of the word $u$. A~letter
$z \in {\cal A}$ is called a~{\em left extension} of a~factor $w$
of the word $u$ if $zw$ belongs to ${\cal L}(u)$. The factor $w$
is a~{\em left special factor} of $u$ if $w$ has more than one
left extensions. Similarly we define right special factors. The
{\em factor complexity} (or simply {\em complexity}) of an
infinite word $u$ is the~function ${\mathcal C}: \mathbb N
\rightarrow \mathbb N$ such that
$$
{\mathcal C}(n)=\mbox{the number of different factors of $u$ of
length $n$}.
$$

We will denote by $\overline{w}$ the reversal of $w$, i.e., if
$w=w_1w_2\cdots w_n$, then $\overline{w}=w_n\cdots w_2w_1$. A~{\em
palindrome} is a~word $w$ which is equal to its reversal
$\overline{w}$. The {\em palindromic complexity} of an infinite
word $u$ is the~function ${\mathcal P}: \mathbb N \rightarrow
\mathbb N$ such that
$$
{\mathcal P}(n)=\mbox{the number of different palindromic factors
of $u$ of length $n$}.
$$
If for a~palindrome $w\in{\cal L}(u)$ there exists $z \in \cal A$
such that $zwz$ is also a factor of $u$, then $z$ is a~{\em
palindromic extension} of $w$. The language of $u$ is {\em closed
under reversal} if for every factor $w$ of $u$, also
$\overline{w}$ is in ${\cal L}(u)$.

A~{\em substitution} on $\cal A^{*}$ is a~morphism $\varphi:{\cal
A^{*}} \rightarrow {\cal A^{*}}$ such that there exists a~letter
$z \in \cal A$ and a~non-empty word $w \in {\cal A}^{*}$
satisfying $\varphi(z)=zw$ and $\varphi(y) \not = \varepsilon$ for
all $y \in \cal A$. Since any morphism satisfies
$\varphi(vw)=\varphi(v)\varphi(w)$ for all $v,w \in {\cal A^{*}}$,
it suffices to define the substitution over the alphabet $\cal A$.
An infinite word $u$ is said to be a~fixed point of the
substitution $\varphi$ if it fulfills
\begin{equation} \label{subst}
u=u_1u_2u_3\cdots =\varphi(u_1)\varphi(u_2)\varphi(u_3)\cdots
=\varphi(u)\,.
\end{equation}
It is obvious that a substitution $\varphi$ has at least one fixed
point, namely $\lim_{n \to \infty}\varphi^{n}(z).$

A~substitution $\varphi$ over the alphabet $\cal A$ is called
primitive if there exists $k \in \mathbb N$ such that for any $z
\in \cal A$ the word $\varphi^{k}(z)$ contains all the letters of
$\cal A$. It has been proved~\cite{DaZa} that if $u$ is a~fixed
point of a~primitive substitution $\varphi$, then $u$ is uniformly
recurrent, that is, for every $n \in \mathbb N$ there exists
$R(n)>0$ such that any factor of $u$ of length $\geq R(n)$
contains all the factors of $u$ of length $n$.

\subsection{Beta-expansions, beta-integers, and R\'enyi expansion of unity} \label{beta}

Let $\beta >1$ be a~real number. For any non-negative real number
$x$ there exist a series of coefficients $(x_i)_{i=-\infty}^k$
such that $x_i \in \mathbb N$  satisfying
$$
0\leq \sum_{i=N}^{k} x_i\beta^{i}< \beta^N \qquad\hbox{ for all }
\ N\in\Z, \ N\leq k\,.
$$
The representation of $x$ in the form of the convergent series
$x=\sum_{i=-\infty}^{k} x_i\beta^{i}$ is called the~{\em
$\beta$-expansions} of $x$. Note that the elements of the sequence
$(x_i)_{i \leq k}$ can also be obtained using the mapping
$T_\beta:[0,1]\mapsto[0,1)$ defined by $T_{\beta}(x):=\beta x -
\lfloor \beta x \rfloor$, by the prescription $x_i = \bigl\lfloor
\beta T_{\beta}^{k-i}(x/ \beta^{k+1}) \bigr\rfloor$. This implies
that the `digits' $x_i$ of the $\beta$-expansion take values in
the set $\{0,1,\dots,\lceil\beta\rceil-1\}$.

Numbers $x\geq 0$ whose $\beta$-expansion is of the form
$\sum_{i=0}^k x_i\beta^{i}$ are called  {\em non-negative
$\beta$-integers}, their set is denoted ${\mathbb Z}_{\beta}^+$.
The set of $\beta$-integers is $\Z_\beta=-{\mathbb
Z}_{\beta}^+\cup {\mathbb Z}_{\beta}^+$. The properties of
$\beta$-integers depend on the so-called {\em R\'enyi expansion of
unity} in base $\beta$ defined as
$$
d_{\beta}(1)=t_1t_2t_3\cdots, \quad \mbox{where }\quad
t_i:=\lfloor \beta T_{\beta}^{i-1}(1)\rfloor\,.
$$
Note that $t_1=\lfloor \beta\rfloor \geq1$. It is interesting to
mention that not every sequence of nonnegative integers is equal
to $d_{\beta}(1)$ for some $\beta$. Parry studied this problem in
his paper~\cite{Parry} and showed that a~sequence $(t_i)_{i \geq
1}$, $t_i \in \mathbb N$, is the R\'enyi expansion of unity for
some number $\beta$ if and only if the sequence satisfies
$$
t_jt_{j+1}t_{j+2}...\prec t_1t_2t_3...\quad \mbox{for every $j
>1$,}
$$
where $\prec$ stands for strictly lexicographically smaller.

\subsection{Infinite words associated with beta-integers}

If $\beta$ is an integer, then $\mathbb Z_{\beta}=\mathbb Z$. For
$\beta \not \in \mathbb N$, Thurston~\cite{Thurston} has shown
that the distances between consecutive $\beta$-integers take
values in the set $\{\Delta_{k}\bigm | k \in \mathbb N \}$, where
\begin{equation}\label{distance}
\Delta_{k}:=\sum_{i=1}^{\infty}\frac{t_{i+k}}{{\beta}^{i}} \
\hbox{for} \ k \in \mathbb N_0.
\end{equation}
It is evident that the set $\{\Delta_{k}\bigm | k \in \mathbb N_0
\}$ is finite if and only if $d_{\beta}(1)$ is eventually
periodic. Numbers $\beta$ with eventually periodic R\'enyi
expansion of unity are called {\em Parry numbers}. Coding the
distances $\Delta_i$ by letters, the set ${\mathbb Z}_{\beta}^{+}$
of non-negative integers can be naturally represented by an
infinite word over a finite alphabet; such a right-sided infinite
word is usually denoted by $u_\beta$.

If the R\'enyi expansion of unity in base $\beta$ is finite
(i.e.,\ $d_\beta(1)$ ends in infinitely many 0's), $\beta$ is said
to be a~{\em simple Parry number}. In this paper, we focus on
non-simple Parry numbers $\beta$, i.e., numbers having eventually
periodic, but infinite R\'enyi expansion of unity. Let $\beta$ be
a~non-simple Parry number and let $m,p$ be minimal such that
\begin{equation}\label{non-simpleRenyi}
d_{\beta}(1)=t_1\dots t_m(t_{m+1}\dots t_{m+p-1})^{\omega}
\end{equation}
is the R\'enyi expansion of unity in base $\beta$. In this case,
there are $m+p-1$ different distances between neighboring
$\beta$-integers. If we assign the letter $i$ to the distance
$\Delta_i$, then ${\mathbb Z}_{\beta}^{+}$ is represented by an
infinite word $u_\beta$ over the alphabet $i \in
\{0,\dots,m+p-1\}$. It has been proved~\cite{Fabre} that the
infinite word $u_\beta$ is the unique fixed point of the following
primitive substitution:
\begin{equation} \label{non-simple}
\begin{array}{rcl} \varphi(j)&= &0^{t_j}(j+1) \  \qquad \hbox{
for }  j\in \{0,1,\dots,
m+p-2\}\,,\\[2mm]
\varphi(m\!+\!p\!-\!1)&=&0^{t_{m+p}}m\,.
\end{array}
\end{equation}

\section{Languages closed under reversal for Parry numbers} \label{reversal}

The existence of palindromes in an infinite word is strongly
related with invariance under reversal of the corresponding
language. This property singles out among all the words $u_\beta$
with non-simple Parry number $\beta$ those, for which $\beta$ is a
quadratic number.

We remind two results; at first, a proposition, proved
in~\cite{AmFrMaPe1}, which explains how the number of palindromes
is connected with the invariance of the language ${\cal
L}(u_{\beta})$ under reversal.

\begin{proposition} \label{reversal_closed}
Let $u$ be an infinite uniformly recurrent word. If for every  $n
\in \mathbb N$ there exists a~palindrome in ${\cal L }(u)$ of
length greater than $n$, then ${\cal L }(u)$ is closed under
reversal.
\end{proposition}

It follows immediately that if there exists a~factor $w \in {\cal
L }(u)$ such that its reversal $\overline{w} \not \in {\cal L
}(u)$, then there is only a~finite number of palindromes in ${\cal
L }(u)$. The following proposition taken from~\cite{Be} determines
the non-simple Parry numbers $\beta$, for which the language
${\cal L}(u_{\beta})$ is closed under reversal.

\begin{proposition}\label{Be}
Let $\beta$ be a non-simple Parry number with  the R\'enyi
expansion of unity $d_\beta(1)=t_1\dots t_m(t_{m+1}\dots
t_{m+p-1})^{\omega}$, where $m,p$ are minimal, and let $u_{\beta}$
be the fixed point of the substitution~(\ref{non-simple}). The
language ${\cal L}(u_{\beta})$ is closed under reversal if and
only if $m=p=1$.
\end{proposition}

The condition $m=p=1$ in the above proposition implies that the
R\'enyi expansion of unity of the number $\beta$ is of the form
\begin{equation}\label{e:quadr}
d_{\beta}(1)=ab^{\omega}\,, \qquad\hbox{ where } \ a-1\geq b \geq
1\,.
\end{equation}
Consequently, $\beta$ is the larger root of the polynomial
$x^2-(a+1)x+a-b$, the substitution $\varphi$ is defined over the
alphabet $\{0,1\}$ by
\begin{equation} \label{subst2}
\varphi(0)=0^a1\,,\quad \varphi(1)=0^b1\,,
\end{equation}
and its fixed point $u_\beta=\lim_{n\to\infty}\varphi^n(0)$ is a
right-sided infinite word of the form
\begin{equation}\label{e:u}
u_\beta = \underbrace{0^a1 \ \cdots \ 0^a1}_{a \hbox{\scriptsize\
times}} \ 0^b1 \ \underbrace{0^a1 \ \cdots \ 0^a1}_{a
\hbox{\scriptsize\ times}} \ 0^b1 \ \cdots
\end{equation}

As a result of Propositions~\ref{reversal_closed} and~\ref{Be}, if
$\beta$ is a non-simple Parry number which is not quadratic, then
the language of the infinite word $u_\beta$ contains only a finite
number of palindromes. Therefore in the study of palindromic
complexity of infinite words $u_\beta$, we will be particularly
interested in quadratic non-simple Parry numbers with $d_\beta(1)$
of the form~\eqref{e:quadr}.

The case of $b=a-1$ has been already studied and the palindromic
complexity of the corresponding infinite word $u_\beta$ is known.
If $b=a-1$, then $\beta$ is the larger root of the polynomial
$x^2-(a+1)x+1$, i.e., $\beta$ is a~quadratic Pisot unit. In
\cite{BuFrGaKr}, it has been shown that the complexity of
$u_\beta$ in this case is ${\mathcal C}(n)=n+1$ for $n\in\N$,
i.e., ${\mathcal C}(n)$ is the smallest possible for aperiodic
words. Infinite words having the smallest possible complexity are
called {\em Sturmian words}. It has been shown in~\cite{Damanik}
that each palindrome $p$ of a~Sturmian word $u$ has exactly one
palindromic extension, i.e.,\ there exists exactly one letter $z
\in \{0,1\}$, such that $zpz \in {\cal L }(u)$. Consequently,
there are exactly 2 palindromes of any odd length and one
palindrome of any even length, i.e.,
\begin{equation}
P(2n)=P(0)=1 \quad \mbox{and} \quad P(2n+1)=P(1)=2 \  \quad
\mbox{for all} \ n \in \mathbb N\,.
\end{equation}

From now on, we will limit our considerations to quadratic
non-simple and non-Sturmian Parry numbers with the R\'enyi
expansion of unity equal to $ d_{\beta}(1)=ab^{\omega}$, where
$a-1>b \geq 1$.

\section{Factor complexity of $u_{\beta}$ } \label{complexityab}

Let us recall some results from~\cite{Ba,FrMaPe2} concerning the
factor complexity of the infinite words $u_\beta$ for a~quadratic
non-simple Parry number $\beta$. It turns out that the notions,
tools, and ideas developed for determining the factor complexity
of $u_\beta$ will be useful also in determining its palindromic
complexity.

Recall that in order to describe the factor complexity of an
infinite word, it suffices to describe its left special factors.
For the case of an infinite word over a~binary alphabet, the
reason is particularly simple. Denoting by $M_n$ the set of all
left special factors of length $n$, it is obvious that the first
difference of complexity satisfies
$$
\Delta C(n)=C(n+1)-C(n)=\#M_n\,.
$$

The first observation about special factors is presented by the
following lemma.

\begin{lemma}\label{l:bloky}
Let $10^r1$ be a~factor of $u_\beta$, then either $r=a$ or $r=b$.
Let $v$ be a~left (right) special factor of $u_\beta$ containing
the letter $1$, then $v$ has the prefix $0^b1$ (the suffix
$10^b$).
\end{lemma}

\begin{proof}
The first statement follows immediately from the form of
substitution~(\ref{subst2}). See~(\ref{e:u}) to understand that
there are no other blocks of the form $10^r1$. Consequently, a
left special factor must have the prefix $0^b1$, otherwise, it
would have a unique left extension.
\end{proof}

Taking a factor $w$ of $u_\beta$, by Lemma~\ref{l:bloky}, the word
$0^b1\varphi(w)0^b$ is also a factor of $u_\beta$. On the other
hand, every factor of $u_\beta$ with the prefix $0^b1$ and suffix
$10^b$ can be written as $0^b1\varphi(w)0^b$ for some factor $w$
of $u_\beta$. Thus, defining a mapping $T:\{0,1\}^*\to\{0,1\}^*$
which will be very useful throughout the paper, we have
\begin{equation}\label{e:T}
w\in{\cal L}(u_\beta) \qquad\Longleftrightarrow\qquad
T(w):=0^b1\varphi(w)0^b\ \in\ {\cal L}(u_\beta)\,.
\end{equation}

For the description of the left special factors of $u_\beta$, one
introduces the following notions.

\begin{definition}
Let $u_{\beta}$ be the infinite word associated with
$d_{\beta}(1)=ab^{\omega}, \ a-1>b \geq 1.$
\begin{itemize}
\item
An infinite word $v$ is called an infinite left special branch of
$u_{\beta}$ if each prefix of $v$ is a~left special factor of
$u_{\beta}$.
\item
A left special factor $w \in {\cal L}(u_{\beta})$ is called
maximal if neither $w0$ nor $w1$ are left special.
\item
A factor $w$ of $u_{\beta}$ is called total bispecial if both $w0$
and $w1$ are left special factors of $u_{\beta}$.
\end{itemize}
\end{definition}

The following facts, proved in~\cite{FrMaPe2}, represent the
consecutive steps needed for the determination of the factor
complexity of $u_\beta$:

\medskip
\begin{trivlist}
\item[\bf -]
Every left special factor is a~prefix of a~maximal left special
factor or of an infinite left special branch. A total bispecial
factor is a common prefix of an infinite left special branch and a
maximal left special factor.

\medskip
\item[\bf -]
The only maximal left special factor which does not contain the
letter 1 has the form $0^{a-1}$ and the only total bispecial
factor containing no letter 1 is $0^{b}$.

\medskip
\item[\bf -]
Every maximal left special (total bispecial) factor of $u_\beta$
that contains at least one letter 1 is also right special, and
therefore by Lemma~\ref{l:bloky}, can be written as
$0^b1\varphi(w)0^b$ for some maximal left special (total
bispecial) factor $w$.

\medskip
\item[\bf -]
If $w \in {\cal L}(u_{\beta})$,  then $w$ is a~left special factor
if and only if $T(w)$ is a~left special factor. Moreover, $w$ is
maximal if and only if $T(w)$ is maximal and $w$ is total
bispecial if and only if $T(w)$ is total bispecial.

\medskip
\item[\bf -]
From the previous, all maximal left special factors have the form:
\begin{equation} \label{U_n}
\begin{array}{rcl}
U^{(1)}&=&0^{a-1}, \\[2mm]
U^{(n)}&=&T(U^{(n-1)}) \ \quad \mbox{for $n \geq 2$.}
\end{array}
\end{equation}

\medskip
\item[\bf -]
All total bispecial factors have the form:
\begin{equation} \label{V_n}
\begin{array}{rcl}
V^{(1)}&=&0^{b}, \\[2mm]
V^{(n)}&=&T(V^{(n-1)}) \ \quad \mbox{for $n \geq 2$.}
\end{array}
\end{equation}

\medskip
\item[\bf -]
$V^{(n-1)}$ is a~prefix of $V^{(n)}$, $V^{(n)}$ is a~prefix of
$U^{(n)}$, and
\begin{equation}\label{e:nerovnost}
|V^{(n)}|<|U^{(n)}|<|V^{(n+1)}| \ \quad \hbox{ for all } \ n \in
\mathbb N\,.
\end{equation}

\medskip
\item[\bf -]
There exists a unique infinite left special branch of $u_\beta$,
namely the infinite word $\lim_{n \to \infty}V^{(n)}$.
\end{trivlist}

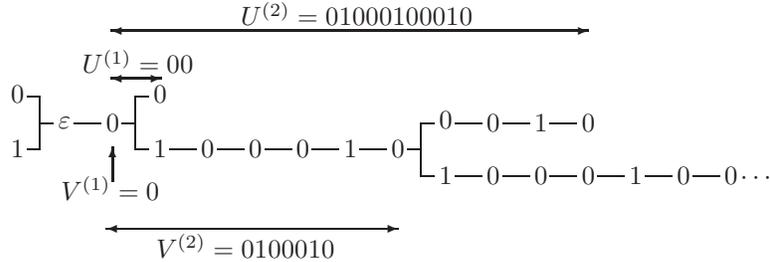
\begin{figure}[ht]
\begin{center}
\begin{picture}(300,110)
\put(26,55){\line(-1,0){5}} \put(21,65){\line(-1,0){5}}
\put(21,45){\line(-1,0){5}} \put(21,45){\line(0,1){20}}
\put(10,63){\mbox{$0$}} \put(10,42){\mbox{$1$}}
\put(28,53){\mbox{$\varepsilon$}} \put(34,55){\line(1,0){10}}
\put(46,52){\mbox{$0$}} \put(52,55){\line(1,0){5}}
\put(57,65){\line(1,0){5}} \put(57,45){\line(1,0){5}}
\put(57,45){\line(0,1){20}}
\put(64,63){\mbox{$0$}} \put(64,42){\mbox{$1$}}
\put(70,45){\line(1,0){10}} \put(82,42){\mbox{$0$}}
\put(88,45){\line(1,0){10}} \put(100,42){\mbox{$0$}}
\put(106,45){\line(1,0){10}} \put(118,42){\mbox{$0$}}
\put(124,45){\line(1,0){10}} \put(136,42){\mbox{$1$}}
\put(142,45){\line(1,0){10}} \put(154,42){\mbox{$0$}}
\put(160,45){\line(1,0){5}} \put(165,55){\line(1,0){5}}
\put(165,35){\line(1,0){5}} \put(165,35){\line(0,1){20}}
\put(172,53){\mbox{$0$}} \put(172,32){\mbox{$1$}}
\put(178,35){\line(1,0){10}} \put(190,32){\mbox{$0$}}
\put(196,35){\line(1,0){10}} \put(208,32){\mbox{$0$}}
\put(214,35){\line(1,0){10}} \put(226,32){\mbox{$0$}}
\put(232,35){\line(1,0){10}} \put(244,32){\mbox{$1$}}
\put(250,35){\line(1,0){10}} \put(262,32){\mbox{$0$}}
\put(268,35){\line(1,0){10}} \put(280,32){\mbox{$0$}}
\put(286,32){\mbox{$\cdots$}}
\put(178,55){\line(1,0){10}} \put(190,52){\mbox{$0$}}
\put(196,55){\line(1,0){10}} \put(208,52){\mbox{$1$}}
\put(214,55){\line(1,0){10}} \put(226,52){\mbox{$0$}}
\put(48.5,33){\vector(0,1){13}} \put(29,26){$V^{(1)}=0$}
\put(53,72){\vector(-1,0){5}} \put(48.5,72){\vector(1,0){18.5}}
\put(37,74){$U^{(1)}=00$}
\put(53,90){\vector(-1,0){5}} \put(48.5,90){\vector(1,0){180}}
\put(97,92){$U^{(2)}=01000100010$}
\put(53,15){\vector(-1,0){7}} \put(48.5,15){\vector(1,0){108}}
\put(65,4){$V^{(2)}=0100010$}
\end{picture}
\end{center}
\caption{Illustration of the tree of left special factors for
$u_{\beta}$ being the fixed point of the substitution
$\varphi(0)=0001, \varphi(1)=01$. We can see total bispecial
factors $V^{(k)}$ and maximal left special factors $U^{(k)}$ for
$k=1,2$.} \label{max_left_complexity}
\end{figure}

If we sum up the previous results, we notice that there is one
left special factor of length $n$ being the prefix of the infinite
left special branch for any $n \in \mathbb N$. Since $|V^{(k)}| <
|U^{(k)}|< |V^{(k+1)}|$ and $V^{(k)}$ is a~prefix of $U^{(k)}$,
there exists a~left special factor of length $n$ being prefix of
$U^{(k)}$ and not prefix of $V^{(k)}$ if $|V^{(k)}|<n \leq
|U^{(k)}|$. Thus, we obtain the following theorem.

\begin{thm}\label{t:faccom}
Let $u_{\beta}$ be the fixed point of the substitution
$\varphi(0)=0^{a}1$, $\varphi(1)=0^{b}1$ for $a-1>b \geq 1$. Then
for all $n \in \mathbb N$,
$$
\Delta C(n)=\left\{
\begin{array}{cl}
2 & \hbox{if }\ |V^{(k)}|<n\leq |U^{(k)}| \ \mbox{ for some~$k \in \mathbb N$},\\
1 &\hbox{otherwise}.\\
\end{array}
\right.
$$
\end{thm}

\section{Palindromic complexity of $u_{\beta}$}
\label{palindromic_complexity}

In this section, we reach the main goal of this paper which is to
determine the palindromic complexity ${\mathcal P}(n)$ of infinite
words associated with quadratic non-simple Parry numbers $\beta$
with the R\'enyi expansion $d_\beta(1)=ab^{\omega}$. As we have
mentioned, the case $b=a-1$ has been already solved, the infinite
word $u_\beta$ being Sturmian under this assumption. Thus we limit
our considerations to the case of $a-1>b$. We proceed in more
steps: First, we investigate palindromic extensions of palindromic
factors of $u_{\beta}$. Second, we study centers of palindromes,
and, finally, we focus on infinite palindromic branches. This
allows us to conclude by determination of the palindromic
complexity.

\subsection{Palindromic extensions of palindromic factors of $u_{\beta}$}

Since the word $u_\beta$ is defined over a binary alphabet, a
palindrome in $u_\beta$ can have two, one, or no palindromic
extensions. It turns out that most usually, palindromes in
$u_\beta$ have exactly one palindromic extension. In this section,
 we determine the exceptional palindromes,
i.e.,\ those that do not have any palindromic extensions (called
maximal palindromes), and those that have two palindromic
extensions.

\begin{definition}
Let $p$ and $q$ be palindromes in ${\cal L}(u_{\beta})$.
\begin{itemize}
\item The palindrome $p$ is a~central factor of the palindrome $q$ if
there exists a~finite word $w \in {\cal L}(u_{\beta})$ such that
$q=wp\overline{w}$.
\item
The palindrome $p$ is called maximal if neither $0p0$ nor $1p1$
belong to ${\cal L}(u_{\beta})$, i.e.,\ if the only palindrome
having $p$ as its central factor is $p$ itself.
\end{itemize}
\end{definition}

Every palindrome is either a~central factor of a~maximal
palindrome or a~central factor of a sequence of palindromes with
increasing length. In the latter case, we will speak about
infinite palindromic branches. The exact definition and
description of infinite palindromic branches in $u_\beta$ follows
in Section~\ref{s:palbr}.

Let us now determine all maximal palindromes in $u_\beta$. First,
realize a simple consequence of Lemma~\ref{l:bloky}.

\begin{lemma}\label{l:rozsireni}
Let $p$ be a palindrome of $u_\beta$ containing at least one
letter 1. Then $p$ is a~central factor of a~palindrome with prefix
$0^b1$ (and, equivalently, suffix $10^b$).
\end{lemma}

Consequently, a palindrome which has not the prefix $0^b1$ is
either not maximal, or does not contain the letter 1. Such
palindrome must be of the form $0^r$ for some $1\leq r\leq a$, and
it can be easily shown that the only maximal palindrome of such
form is $0^{a-1}$.


In order to determine all maximal palindromes of $u_\beta$, let us
now study the behaviour of palindromes with respect to the
substitution $\varphi$. Since every maximal palindrome $p$
containing the letter 1 has the prefix $0^b1$ and the suffix
$10^b$, one can always write it as $0^b1\varphi(q)0^b$ for some
factor $q\in{\cal L}(u_\beta)$. The following lemma shows that
similarly as in the case of left special factors studied for
determining the factor complexity, we can use the substitution
$\varphi$ to reduce all maximal palindromes containing at least
one letter 1 to the maximal palindrome $0^{a-1}$.


\begin{lemma} \label{pT}
Let $p \in {\cal L}(u_{\beta})$. Then $p$ is a~palindrome in
$u_\beta$ if and only if $T(p)=0^b 1\varphi(p)0^b$ is a~palindrome
in $u_\beta$. Moreover,
$$
\bigl\{z\in \{0,1\} \bigm | zpz \in {\cal L}(u_{\beta})\bigr\}
=\bigl\{z\in \{0,1\} \bigm | zT(p)z \in {\cal
L}(u_{\beta})\bigr\}\,.
$$
\end{lemma}

\begin{proof}
Let us investigate under which condition $0^b 1\varphi(p)0^b$ is a
palindrome in $u_{\beta}$.  We will study for which $p$ it holds
that
$$
0^b\overline{\varphi(p)}10^b=0^b 1\varphi(p)0^b\,.
$$
Note that we can write
$\overline{\varphi(0)}=10^a=1\varphi(0)1^{-1}$ and
$\overline{\varphi(1)}=10^b=1\varphi(1)1^{-1}$. Then denoting
$p=p_1...p_k$, we have
$$
\overline{\varphi(p)}1=\overline{\varphi(p_1)\dots \varphi(p_k)}1=
\overline{\varphi(p_k)}\dots
\overline{\varphi(p_1)}1=1\varphi(\overline{p})\,.
$$
As $\varphi(p)=\varphi(\overline{p})$ if and only if
$p=\overline{p}$, the first statement of the lemma is proved.

\medskip
Let us now verify $\{z \in \{0,1\} \bigm | zpz \in {\cal
L}(u_{\beta})\} =\{z\in \{0,1\} \bigm | zT(p)z \in {\cal
L}(u_{\beta})\}$.
\begin{itemize}
\item[$\subset$]
Let $0p0 \in {\cal L}(u_{\beta})$. Then
$T(0p0)=0^b1\varphi(0)\varphi(p)\varphi(0)0^b \in {\cal
L}(u_{\beta})$. Since $\varphi(0)$ has the suffix $00^{b}1$ and
the prefix $0^{b}0$, we have $0T(p)0 \in {\cal L}(u_{\beta})$.

Let now $1p1 \in {\cal L}(u_{\beta})$. Then
$T(1p1)=0^b1\varphi(1)\varphi(p)\varphi(1)0^b \in {\cal
L}(u_{\beta})$ and from $\varphi(1)=0^b1$,  we deduce that $1T(p)1
\in {\cal L}(u_{\beta})$.

\item[$\supset$]
Let $0T(p)0\in {\cal L}(u_{\beta})$. Due to Lemma~\ref{l:bloky},
$0T(p)0=0^{b+1}1\varphi(p)0^{b+1}$ can be uniquely extended on
both sides to the factor
$$
0^b10^a1\varphi(p)0^a10^b=0^b1\varphi(0p0)0^b=T(0p0) \in {\cal
L}(u_{\beta})\,,
$$
which implies $0p0 \in {\cal L}(u_{\beta})$.

Similarly, let
$1T(p)1\in {\cal L}(u_{\beta})$. Then 
$1T(p)1=10^b1\varphi(p)0^b1$ is a central factor of
$0^b10^b1\varphi(p)0^b10^b = 0^b1\varphi(1p1)0^b=T(1p1) \in {\cal
L}(u_{\beta})$, which implies $1p1\in {\cal L}(u_{\beta})$.
\end{itemize}
\end{proof}

Using Lemma~\ref{pT}, we can see that the maximal left special
factors $U^{(n)}$, defined by~\eqref{U_n}, are at the same time
maximal palindromes.

\begin{proposition}
Let $p$ be a palindrome in $u_\beta$. Then $p$ is a~maximal
palindrome if and only if $p=U^{(n)}$ for a~positive integer $n$.
\end{proposition}

\pf From their definition by~\eqref{U_n} and from Lemma~\ref{pT},
every $U^{(n)}$ is a maximal palindrome. We show the opposite
implication by induction on the number of letters 1 contained in
$p$. If $p$ does not contain 1 and $p$ is maximal, then $p$ is
equal to $U^{(1)}= 0^{a-1}.$ There is no maximal palindrome
containing only one letter 1. Let $p$ contain at least two letters
1. Then by Lemma~\ref{l:rozsireni}, $p$ is of the form
$p=0^b1\varphi(q)0^b$. Using Lemma~\ref{pT}, we know that $q$ is
a~maximal palindrome and $q$ contains a smaller number of letters
1 than $p$. By induction hypothesis, $q=U^{(k)}$ for some $k \in
\mathbb N$ and $p=0^b1\varphi(U^{(k)})0^b=U^{(k+1)}$. \pfk

\bigskip
Let us now determine the palindromes which have two palindromic
extensions. Such a palindrome is both a left special and a right
special factor of $u_\beta$. Therefore by Lemma~\ref{l:bloky}, it
must have the prefix $0^b1$ and the suffix $10^b$, whenever it
contains at least one letter 1. Thus again, we can write such a
palindrome $p$ as $p=0^b1\varphi(q)0^b$ for some $q\in{\cal
L}(u_\beta)$. Now it suffices to realize that the only palindrome
of the form $0^r$ with two palindromic extensions, is
$V^{(1)}=0^b$. Using Lemma~\ref{pT}, we deduce that the total
bispecial factors $V^{(n)}$ defined by~\eqref{V_n} have two
palindromic extensions.

\begin{proposition}
Let $p$ be a palindrome in $u_\beta$. Then $p$ has two palindromic
extensions if and only if $p=V^{(n)}$ for a~positive integer $n$.
\end{proposition}

\pf From Lemma~\ref{pT}, it is clear that palindromes $V^{(n)}$
have two palindromic extensions. It suffices to show that any
palindrome with two palindromic extensions is equal to some
$V^{(n)}$. Let $p$ be a~palindrome such that $0p0, 1p1 \in
Pal(u_{\beta})$, then $p$ is a~left special factor and it is not
maximal. Without loss of generality, we can suppose that $p0$ is
a~left special factor, i.e., $0p0, 1p0 \in {\cal L}(u_{\beta})$ .
Using the fact that ${\cal L}(u_{\beta})$ is closed under
reversal, we have $0p1, 1p1 \in {\cal L}(u_{\beta})$. It means
that also $p1$ is a~left special factor. Consequently, $p$ is
a~total bispecial factor, i.e., $p=V^{(k)}$ for some~$k \in
\mathbb N$. \pfk

Let us summarize the results of this section by classifying
palindromes in $u_\beta$ according to their number of palindromic
extensions. This a crucial step for the determining of palindromic
complexity of $u_\beta$. It is interesting to notice that
sequences $U^{(n)}$ of maximal left special factors and $V^{(n)}$
of total bispecial factors play an mportant role not only in
computing the factor complexity, but also in determining the
palindromic complexity. For illustration, consider the example
from Figure~\ref{max_left_complexity} in
Section~\ref{complexityab} and check that the total bispecial
factors and maximal left special factors illustrated there are
indeed palindromes.

\begin{proposition} \label{Pal}
Let $p$ be a palindrome in $u_{\beta}$. Then
\begin{enumerate}
\item
$p$ is a~maximal palindrome $\Leftrightarrow$ $p=U^{(n)}$ for
a~positive integer $n$.
\item
$p$ has two palindromic extensions $\Leftrightarrow$ $p=V^{(n)}$
for a~positive integer $n$.
\item
$p$ has one palindromic extension $\Leftrightarrow$ $p \not
=U^{(n)} \wedge p\not=V^{(n)}$ for all $n \in \mathbb N$.
\end{enumerate}
\end{proposition}

Using the inequality $|V^{(n)}| < |U^{(n)}| < |V^{(n+1)}|$ for all
$n \in \mathbb N$, we obtain the following corollary of
Proposition~\ref{Pal}.

\begin{corollary}
Let $u_{\beta}$ be the fixed point of the substitution
$\varphi(0)=0^{a}1, \ \varphi(1)=0^{b}1, \ a-1>b\geq 1.$ Then
$$
{\mathcal P}(n+2)-{\mathcal P}(n)=\left\{
\begin{array}{rl}
1 & \hbox{ if }\ n=|V^{(k)}| \ \mbox{for a~positive integer $k$},\\
-1& \hbox{ if }\ n=|U^{(k)}| \ \mbox{for a~positive integer $k$},\\
0 &\mbox{otherwise}.\\
\end{array}
\right.
$$
\end{corollary}


Combining with Theorem~\ref{t:faccom}, we can derive a simple
connection of palindromic complexity with the second difference of
factor complexity:
$$
\Delta^{2}{\cal C}(n)=\Delta {\cal C}(n+1)-\Delta {\cal
C}(n)={\cal P}(n+2)-{\cal P}(n)\,.
$$
This relation is essential to show the following corollary.

\begin{corollary}\label{c:rovnost}
Let $u_{\beta}$ be the fixed point of the substitution
$\varphi(0)=0^{a}1, \ \varphi(1)=0^{b}1, \ a-1>b\geq1.$ Then
\begin{equation}\label{e:vztah}
{\cal P}(n+1)+{\cal P}(n)=\Delta {\cal C}(n)+2 \qquad \mbox{for
all} \ n \in \mathbb N\,.
\end{equation}
\end{corollary}

\begin{proof}
We have ${\cal P}(n+1)+{\cal P}(n)={\cal P}(0)+{\cal P}(1)
+\displaystyle {\sum_{i=1}^{n}}\bigl ({\cal P}(i+1)-{\cal
P}(i-1)\bigr )=$ $=1+2+ \displaystyle {\sum_{i=1}^{n}}\bigl
(\Delta {\cal C}(i)-\Delta {\cal C}(i-1)\bigr)= 3+\Delta {\cal
C}(n)-{\cal C}(1)+{\cal C}(0)=\Delta {\cal C}(n)+2$.
\end{proof}


Let us mention that for uniformly recurrent infinite words in
general, one has ${\cal P}(n+1)+{\cal P}(n)\leq\Delta {\cal
C}(n)+2$ for all $n\in\N$. This is shown in~\cite{ThCompSci}
together with several examples of infinite words where the
equality is reached. By Corollary~\ref{c:rovnost}, we have found
yet another class of infinite words with this property.

 A~direct
consequence of the above corollary is the fact that the
palindromic complexity of $u_\beta$ is bounded since $P(n+1)+P(n)
\leq 4$.

\subsection{Centers of palindromes}

We have seen that the set of palindromes of $u_{\beta}$ is closed
under the mapping $p\, \mapsto T(p)=0^b1\varphi(p)0^b$. Let us
study how $T$ acts on the central factors of palindromes.

\begin{definition}
Let $p \in {\cal L}(u_{\beta})$ be a~palindrome of odd length. The
center of $p$ is a~letter $z \in \{0,1\}$ such that
$p=wz\overline{w}$, where $w \in {\cal L}(u_{\beta})$. The center
of a~palindrome of even length is the empty word.
\end{definition}

\begin{lemma}\label{centrapTp}
Let $p,q$ be palindromes in $u_\beta$. If $q$ is a~central factor
of $p$, then $T(q)$ is a~central factor of $T(p)$.
\end{lemma}

\begin{proof}
Let $p=w0q0\overline{w}$ for a $w\in{\cal L}(u_\beta)$. Then
$T(p)=0^b1\varphi(w)0^a1\varphi(q)0^a1\varphi(\overline{w})0^b$.
According to Lemma~\ref{pT}, $T(p)$ is a palindrome, and clearly,
$T(q)=0^b1\varphi(q)0^b$ is its central factor. The proof is
similar for $p=w1q1\overline{w}$.
\end{proof}

Now, using Lemma~\ref{centrapTp}, we can describe how the center
of the palindrome $T(p)$ depends on the center of the palindrome
$p$.

\begin{lemma}\label{centra}
Let $p$ be a palindrome in $u_{\beta}$.
\begin{enumerate}
\item[(i)]
If $p$ has the center $\varepsilon$, then $T(p)$ has the center 1.
\item[(ii)]
If $p$ has the center 0, then $T(p)$ has the center
$\left\{\begin{array}{cl}
0 & \mbox{for} \ a \ \mbox{odd},\\
\varepsilon & \mbox{for} \ a \ \mbox{even}.
\end{array}
\right. $
\item[(iii)]
If $p$ has the center 1, then $T(p)$ has the center
$\left\{\begin{array}{cl}
0 & \mbox{for} \ b \ \mbox{odd},\\
\varepsilon & \mbox{for} \ b \ \mbox{even}.
\end{array}
\right. $
\end{enumerate}
\end{lemma}

\begin{proof}
Let us verify for example the statement (ii). Using
Lemma~\ref{centrapTp}, it is evident that if $p$ has the center
$0$, then $T(p)$ has the central factor $T(0)=0^b10^a10^b$.
Consequently, the center of $T(p)$ is either $0$ if $a$ is odd, or
$\varepsilon$ if $a$ is even. Statements (i) and (iii) can be
proved analogously.
\end{proof}

Lemmas~\ref{centrapTp} and~\ref{centra} allow us to describe the
centers of palindromes with two palindromic extensions $V^{(n)}$
and the centers of maximal palindromes $U^{(n)}$.

\begin{proposition} \label{formVn} ~
The centers and central factors of palindromes $V^{(n)}$ with two
palindromic extensions depend on the values of parameters $a,b$.
\begin{enumerate}
\item[(i)]
Let $b$ be even. Then for all $n \in \mathbb N$, $V^{(n)}$ is
a~central factor of $V^{(n+2)}$. Moreover, $V^{(2n)}$ has the
center 1 and $V^{(2n-1)}$ has the center $\varepsilon$.
\item[(ii)]
Let $b$ be odd and $a$ even. Then for all $n \in \mathbb N$,
$V^{(n)}$ is a~central factor of $V^{(n+3)}$. Moreover, $V^{(3n)}$
has the center 1, $V^{(3n-1)}$ has the center $\varepsilon$, and
$V^{(3n-2)}$ has the center 0.
\item[(iii)]
Let both $b$ and $a$ be odd. Then for all $n \in \mathbb N$,
$V^{(n)}$ is a~central factor of $V^{(n+1)}$ and has the center 0.
\end{enumerate}
\end{proposition}

\begin{proof}
In order to show the statement (i), it suffices to verify that
$V^{(1)}$ is a~central factor of $V^{(3)}$ and that $V^{(1)}$ has
the center $\varepsilon$ and $V^{(2)}$ has the center 1. The
statement (i) then follows by induction on $n \in \mathbb N$.
Since $b$ is even, $V^{(1)}=0^b$ has the center $\varepsilon$. By
Lemma~\ref{centra}, $V^{(2)}$ has the center 1. Applying
Lemma~\ref{centrapTp}, one can see that $V^{(3)}$ has the central
factor $T(1)=0^b10^b10^b$, i.e.,\ it has also $V^{(1)}$ as its
central factor.

Proofs of statements (ii) and (iii) are analogous.
\end{proof}

\begin{proposition} \label{formUn} ~
The centers and central factors of maximal palindromes $U^{(n)}$
depend on the values of $a$ and $b$. Let $n \in \mathbb N$, we
have the following cases:
\begin{enumerate}
\item[(i)]
Let $b$ be even and $a$ odd. Then for all $n \in \mathbb N$,
$V^{(n)}$ is a central factor of $U^{(n)}$. Moreover, $U^{(2n-1)}$
has the center $\varepsilon $ and $ U^{(2n)}$ has the center 1.

\item[(ii)]
Let both $b$ and $a$ be even. Then $U^{(1)}=0^{a-1}$ is the only
maximal palindrome with the center 0. For all $n \in \mathbb N$,
$V^{(n)}$ is a central factor of $U^{(n+1)}$. Moreover, $U^{(2n)}$
has the center $\varepsilon $ and $U^{(2n+1)}$ has the center 1.

\item[(iii)]
Let $b$ be odd and $a$ even. Then for all $n \in \mathbb N$,
$V^{(n)}$ is a central factor of $U^{(n)}$. Moreover, $U^{(3n-2)}$
has the center 0, $ U^{(3n-1)}$ has the center $\varepsilon $, and
$U^{(3n)}$ has the center 1.

\item[(iv)]
Let both $b$ and $a$ be odd. The only maximal palindrome with the
center $\varepsilon $ is $U^{(1)}=0^{a-1}$. The only maximal
palindrome having the center 1 is $U^{(2)}$. For $n \geq 3$,
$U^{(n)}$ has the center 0 and the central factor $V^{(n-2)}$.
\end{enumerate}
\end{proposition}

\begin{proof}
Let us show for example the statement (ii). Since $a$ is even,
$U^{(1)}=0^{a-1}$ has the center 0, Lemma~\ref{centra} implies
that $U^{(2)}$ has the center $\varepsilon $ and the central
factor $T(0)=0^b10^a10^b$. Consequently, $V^{(1)}=0^b$ is also
a~central factor of $U^{(2)}$, $b$ being even. Applying
Lemma~\ref{centra}, we obtain that $U^{(3)}$ has the center 1. The
statement follows by induction on $n \in \mathbb N$.
\end{proof}

\subsection{Infinite palindromic branches}\label{s:palbr}

Every palindrome is either a~central factor of a~maximal
palindrome $U^{(n)}$, or a~central factor of an infinite
palindromic branch. Their knowledge is essential for
determination of the palindromic complexity of $u_\beta$.

\begin{definition}
Let $v=v_1v_2v_3\cdots$ be a~right-sided infinite word over the
alphabet $\{0,1\}$. Denote by $\overline{v}$ the left-sided
infinite word $\overline{v}=\cdots v_3v_2v_1$. Let $z \in
\{\varepsilon,0,1\}$. If for every $n \in \mathbb N$, the
palindrome $p=v_nv_{n-1}...v_1zv_1v_2...v_n$ belongs to ${\cal
L}(u_{\beta})$, then the bidirectional infinite word
$\overline{v}zv$ is called an infinite palindromic branch of
$u_{\beta}$ with the center $z$, and the palindrome $p$ is called
the central factor of the infinite palindromic branch
$\overline{v}zv$.
\end{definition}

\begin{lemma}\label{AtMostOneInf}
Let $z \in \{\varepsilon,0,1\}$. Then there exists at most one
infinite palindromic branch with the center $z$.
\end{lemma}

\begin{proof}
It follows from Lemma~\ref{centrapTp} that applying the
substitution $\varphi$ on an infinite palindromic branch, one
obtains again an infinite palindromic branch. More precisely, from
$\overline{v}zv$ one obtains
$\overline{\varphi(v)}10^a1\varphi(v)$ if $z=0$,
$\overline{\varphi(v)}10^b1\varphi(v)$ if $z=1$, and
$\overline{\varphi(v)}1\varphi(v)$ if $z=\varepsilon$.

We show the statement of the lemma by contradiction. Suppose that
there are two different infinite palindromic branches with the
same center $z$. By applying the substitution $\varphi$, we obtain
two different infinite palindromic branches with longer common
central factor. Repeating this procedure, we can construct
infinitely many different infinite palindromic branches, which is
a~contradiction with boundedness of the palindromic complexity.
\end{proof}

We now describe the infinite palindromic branch of $u_\beta$ with
the center $z\in\{\varepsilon,0,1\}$ and the common central
factors of this branch with maximal palindromes having the same
center $z$. With this in hand, we will be able to summarize the
values of the palindromic complexity. Note that the candidate for
the longest common prefix of a maximal palindrome and an infinite
palindromic branch with the same center is a palindrome, which has
two palindromic extensions, thus one of the palindromes $V^{(n)}$.
Using Proposition~\ref{formVn}, we can describe all infinite
palindromic branches of $u_{\beta}$.

\begin{proposition} \label{allbidirect}~
\begin{enumerate}
\item[(i)]
Let $b$ be even and $a$ odd. There exists an infinite palindromic
branch with the center $z$ for all $z\in\{\varepsilon, 1,0\}$,
namely the bidirectional limit:
$$
\begin{array}{ll}
\qquad\qquad \displaystyle{\lim_{n\to\infty}V^{(2n-1)}} & \mbox{having the center $\varepsilon$,}\\
\qquad\qquad \displaystyle{\lim_{n\to\infty}V^{(2n)}} & \mbox{having the center 1,}\\
\qquad\qquad \displaystyle{\lim_{n\to\infty}W^{(n)}} & \mbox{where $W^{(1)}=0, W^{(n)}=T(W^{(n-1)})$,}\\
& \mbox{having the center 0.}
\end{array}
$$
For $n\in\N$, the longest common central factor of a maximal
palindrome $U^{(n)}$ and the infinite palindromic branch with the
same center is $V^{(n)}$. The limit of $W^{(n)}$ has a common
central factor with no maximal palindrome.

\item[(ii)]
Let both $b$ and $a$ be even. There exists an infinite palindromic
branch with the center $z$ for $z\in\{\varepsilon,1\}$, namely the
bidirectional limit:
$$
\begin{array}{ll}
\displaystyle{\lim_{n\to\infty}V^{(2n-1)}} &\ \mbox{having the center $\varepsilon $,}\\
\displaystyle{\lim_{n\to\infty}V^{(2n)}} &\ \mbox{having the
center 1.}
\end{array}
$$
There is no infinite palindromic branch with the center 0. For
$n\in\N$, $n\geq 2$, the longest common central factor of a
maximal palindrome $U^{(n)}$ and the infinite palindromic branch
with the same center is $V^{(n-1)}$.

\item[(iii)]
Let $b$ be odd and $a$ even. There exists an infinite palindromic
branch with the center $z$ for $z\in\{0,\varepsilon,1\}$, namely
the bidirectional limit:
$$
\begin{array}{ll}
\displaystyle{\lim_{n\to\infty}V^{(3n-2)}} &\ \mbox{having the center 0,}\\
\displaystyle{\lim_{n\to\infty}V^{(3n-1)}} &\ \mbox{having the center $\varepsilon $,}\\
\displaystyle{\lim_{n\to\infty}V^{(3n)}} &\ \mbox{having the
center 1.}
\end{array}
$$
For $n\in\N$, the longest common central factor of a maximal
palindrome $U^{(n)}$ and the infinite palindromic branch with the
same center is $V^{(n)}$.

\item[(iv)]
Let both $b$ and $a$ be odd. There exists an infinite palindromic
branch with the center 0, namely the bidirectional limit of
palindromes $V^{(n)}$, $n\in\N$. There is neither an infinite
palindromic branch with the center $\varepsilon$, nor with the
center 1. For $n\in\N$, $n\geq 3$, the longest common central
factor of a maximal palindrome $U^{(n)}$ and the infinite
palindromic branch is $V^{(n-2)}$.
\end{enumerate}
\end{proposition}

\begin{proof}
It follows from Lemma~\ref{AtMostOneInf} that there is at most one
infinite palindromic branch with the center $z$ for each $z \in
\{\varepsilon, 0, 1\}$. In order to verify that the bidirectional
limits of palindromes exist, it suffices to use
Proposition~\ref{formVn}, or Lemma~\ref{centrapTp}, to see that
the palindromes are central factors of one another. Let us explain
why in cases (ii) and (iv) one does not have an infinite
palindromic branch with every center.

(ii) Let $b$ and $a$ be even, suppose that there is a~palindromic
branch with the center 0. Necessarily, this branch has a~block of
the form $0^a$ or $0^b$ as its central factor. It is impossible
owing to the fact that both $a$ and $b$ are even.

(iv) Let both $b$ and $a$ be odd, suppose that there is an
infinite palindromic branch with the center $\varepsilon$. Then it
has a~block of the form $0^a$ or $0^b$ as its central factor. It
is impossible since both $a$ and $b$ are odd. Suppose now that
there exists an infinite palindromic branch with the center 1.
Take a central factor of this palindromic branch of the form
$T(p)$ for a palindrome $p$ containing at least two letters 1.
Using Lemma~\ref{centra}, $p$ must have the center $\varepsilon$,
and thus have a central factor $0^a$ or $0^b$, which is
impossible.

The statements about the maximal common central factor of maximal
palindromes and infinite palindromic branches are a~consequence of
Proposition~\ref{formUn}.
\end{proof}

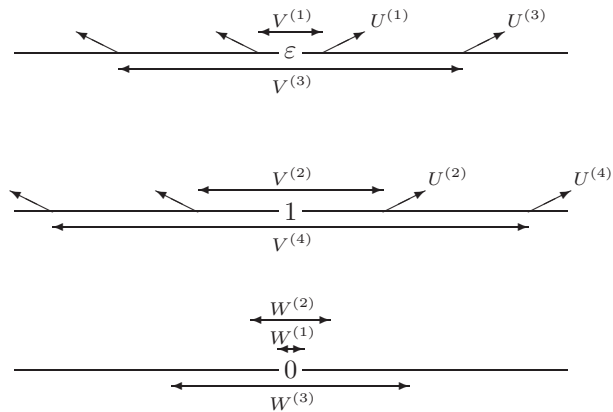
\begin{figure}[!ht]
\begin{center}
\begin{picture}(250,175)
\put(15.5,30){\line(1,0){100}} \put(117.5,27){0}
\put(124.5,30){\line(1,0){100}}
\put(115,38){\vector(1,0){10}} \put(118,38){\vector(-1,0){3}}
\put(112,40){\scriptsize $W^{(1)}$}
\put(105,49){\vector(1,0){30}} \put(108,49){\vector(-1,0){3}}
\put(112,51){\scriptsize $W^{(2)}$}
\put(75,24){\vector(1,0){90}} \put(78,24){\vector(-1,0){3}}
\put(112,15){\scriptsize $W^{(3)}$}
\put(15.5,90){\line(1,0){100}} \put(117.5,87){1}
\put(124.5,90){\line(1,0){100}}
\put(85,98){\vector(1,0){70}} \put(88,98){\vector(-1,0){3}}
\put(113,100){\scriptsize $V^{(2)}$}
\put(30,84){\vector(1,0){180}} \put(33,84){\vector(-1,0){3}}
\put(113,75){\scriptsize $V^{(4)}$}
\put(85,90){\vector(-2,1){16}} \put(155,90){\vector(2,1){16}}
\put(172,100){\scriptsize $U^{(2)}$}
\put(30,90){\vector(-2,1){16}} \put(210,90){\vector(2,1){16}}
\put(227,100){\scriptsize $U^{(4)}$}
\put(15.5,150){\line(1,0){100}} \put(117.5,148){$\varepsilon$}
\put(124.5,150){\line(1,0){100}}
\put(108,158){\vector(1,0){24}} \put(111,158){\vector(-1,0){3}}
\put(113,160){\scriptsize $V^{(1)}$}
\put(55,144){\vector(1,0){130}} \put(58,144){\vector(-1,0){3}}
\put(113,135){\scriptsize $V^{(3)}$}
\put(55,150){\vector(-2,1){16}} \put(185,150){\vector(2,1){16}}
\put(202,160){\scriptsize $U^{(3)}$}
\put(108,150){\vector(-2,1){16}} \put(132,150){\vector(2,1){16}}
\put(150,160){\scriptsize $U^{(1)}$}
\end{picture}
\caption{Illustration of maximal palindromes and infinite
palindromic branches for $b$ even and $a$ odd. There is one
infinite branch with center $\varepsilon$, one with center 1, and
one with center 0. There are infinitely many maximal palindromes
with center $\varepsilon $ and 1.} \label{maximal1}
\end{center}
\end{figure}

\begin{figure}[!ht]
\begin{center}
\begin{picture}(250,140)
\put(115.5,10){\vector(-1,0){20}} \put(117.5,7){0}
\put(124.5,10){\vector(1,0){20}} \put(113,16){\scriptsize
$U^{(1)}$}
\put(15.5,60){\line(1,0){100}} \put(117.5,57){1}
\put(124.5,60){\line(1,0){100}}
\put(85,68){\vector(1,0){70}} \put(88,68){\vector(-1,0){3}}
\put(113,70){\scriptsize $V^{(2)}$}
\put(30,54){\vector(1,0){180}} \put(33,54){\vector(-1,0){3}}
\put(113,45){\scriptsize $V^{(4)}$}
\put(85,60){\vector(-4,1){44}} \put(155,60){\vector(4,1){44}}
\put(202,70){\scriptsize $U^{(3)}$}
\put(15.5,120){\line(1,0){100}} \put(117.5,118){$\varepsilon$}
\put(124.5,120){\line(1,0){100}}
\put(108,128){\vector(1,0){24}} \put(111,128){\vector(-1,0){3}}
\put(113,130){\scriptsize $V^{(1)}$}
\put(55,114){\vector(1,0){130}} \put(58,114){\vector(-1,0){3}}
\put(113,105){\scriptsize $V^{(3)}$}
\put(108,120){\vector(-4,1){35}} \put(132,120){\vector(4,1){35}}
\put(169,130){\scriptsize $U^{(2)}$}
\put(55,120){\vector(-4,1){45}} \put(185,120){\vector(4,1){45}}
\put(232,130){\scriptsize $U^{(4)}$}
\end{picture}

\caption{Illustration of maximal palindromes and infinite
palindromic branches for $a$ and $b$ even. There is one infinite
branch with center $\varepsilon$ and one with center 1. There are
infinitely many maximal palindromes with center $\varepsilon$ and
1. There is only one maximal palindrome with center 0.}
\label{maximal2}
\end{center}
\end{figure}
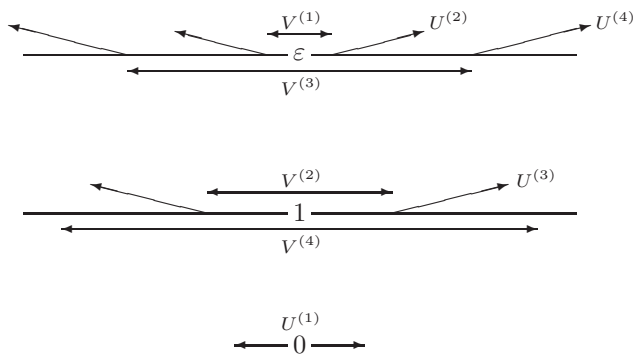

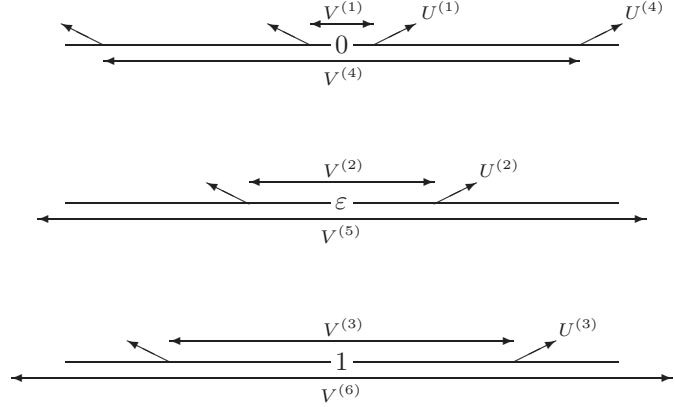
\begin{figure}[!ht]
\begin{center}

\begin{picture}(250,175)
\put(15.5,30){\line(1,0){100}} \put(117.5,27){1}
\put(124.5,30){\line(1,0){100}}
%
%
%
\put(-5,24){\vector(1,0){250}} \put(-2,24){\vector(-1,0){3}}
\put(112,15){\scriptsize $V^{(6)}$}
\put(15.5,90){\line(1,0){100}} \put(117.5,88){$\varepsilon$}
\put(124.5,90){\line(1,0){100}}
\put(85,98){\vector(1,0){70}} \put(88,98){\vector(-1,0){3}}
\put(113,100){\scriptsize $V^{(2)}$}
\put(30,144){\vector(1,0){180}} \put(33,144){\vector(-1,0){3}}
\put(113,135){\scriptsize $V^{(4)}$}
\put(85,90){\vector(-2,1){16}} \put(155,90){\vector(2,1){16}}
\put(172,100){\scriptsize $U^{(2)}$}
\put(30,150){\vector(-2,1){16}} \put(210,150){\vector(2,1){16}}
\put(227,160){\scriptsize $U^{(4)}$}
\put(5,84){\vector(1,0){230}} \put(8,84){\vector(-1,0){3}}
\put(112,75){\scriptsize $V^{(5)}$}
\put(15.5,150){\line(1,0){100}} \put(117.5,147){0}
\put(124.5,150){\line(1,0){100}}
\put(108,158){\vector(1,0){24}} \put(111,158){\vector(-1,0){3}}
\put(113,160){\scriptsize $V^{(1)}$}
\put(55,38){\vector(1,0){130}} \put(58,38){\vector(-1,0){3}}
\put(113,40){\scriptsize $V^{(3)}$}
\put(55,30){\vector(-2,1){16}} \put(185,30){\vector(2,1){16}}
\put(202,40){\scriptsize $U^{(3)}$}
\put(108,150){\vector(-2,1){16}} \put(132,150){\vector(2,1){16}}
\put(150,160){\scriptsize $U^{(1)}$}
\end{picture}

\caption{Illustration of maximal palindromes and infinite
palindromic branches for $b$ odd and $a$ even. There are infinite
branches and infinitely many maximal palindromes with centers
$0,1, \varepsilon$.} \label{maximal3}
\end{center}
\end{figure}

\begin{figure}[!ht]
\begin{center}

\begin{picture}(250,140)
\put(115.5,10){\vector(-1,0){20}} \put(117.5,8){$\varepsilon$}
\put(124.5,10){\vector(1,0){20}} \put(113,16){\scriptsize
$U^{(1)}$}
\put(115.5,50){\vector(-1,0){45}} \put(117.5,47){1}
\put(124.5,50){\vector(1,0){45}}
\put(113,56){\scriptsize $U^{(2)}$}
\put(85,104){\vector(1,0){70}} \put(88,104){\vector(-1,0){3}}
\put(113,95){\scriptsize $V^{(2)}$}
\put(35,90){\vector(1,0){170}} \put(38,90){\vector(-1,0){3}}
\put(113,81){\scriptsize $V^{(4)}$}
\put(85,110){\vector(-4,-1){66}} \put(155,110){\vector(4,-1){66}}
\put(198,127){\scriptsize $U^{(3)}$}
\put(15.5,110){\line(1,0){100}} \put(117.5,107){0}
\put(124.5,110){\line(1,0){100}}
\put(108,118){\vector(1,0){24}} \put(111,118){\vector(-1,0){3}}
\put(113,120){\scriptsize $V^{(1)}$}
\put(55,130){\vector(1,0){130}} \put(58,130){\vector(-1,0){3}}
\put(113,132){\scriptsize $V^{(3)}$}
\put(108,110){\vector(-4,1){65}} \put(132,110){\vector(4,1){65}}
\put(55,110){\vector(-4,1){60}} \put(185,110){\vector(4,1){60}}
\put(222,95){\scriptsize $U^{(4)}$}
\put(246,126){\scriptsize $U^{(5)}$}
\end{picture}

\caption{Illustration of maximal palindromes and infinite
palindromic branches for $a$ and $b$ odd. The only infinite
palindromic branch has center 0. There are infinitely many maximal
palindromes with center 0. There is only one maximal palindrome
with center $\varepsilon$ and one with center 1.} \label{maximal4}
\end{center}
\end{figure}
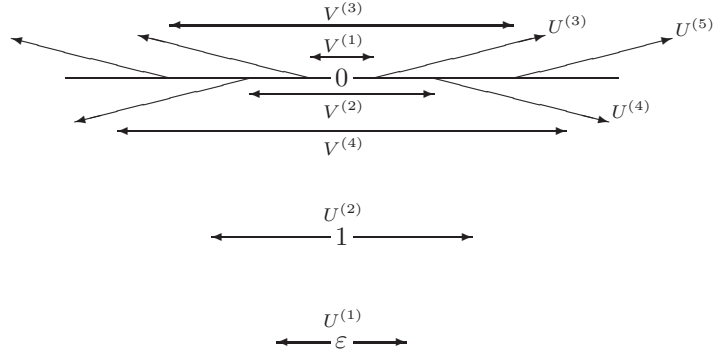

\subsection{Explicit values of the palindromic complexity of $u_{\beta}$}

We are now in position to derive explicitly the values of the
palindromic complexity of the infinite word $u_\beta$ dependingly
on the parity of parameters $a,b$ of the R\'enyi expansion of
unity $d_\beta(1)=ab^\omega$. We have investigated maximal
palindromes, infinite palindromic branches, and their centers.
Determining the complexity is easy with the use of
Figures~\ref{maximal1}--\ref{maximal4}, which visualize the
structure of maximal palindromes and of infinite palindromic
branches, according to Proposition~\ref{allbidirect}.

\begin{thm}
Let $u_{\beta}$ be the fixed point of the substitution
$\varphi(0)=0^{a}1, \ \varphi(1)=0^{b}1, \ a-1>b \geq 1.$ Then 4
cases can appear according to the values of parameters $a$ and
$b$, $n \in \mathbb N \cup \{0\}$:

\begin{enumerate}
\item[(i)]
Let $b$ be even and $a$ odd.
$$
\begin{array}{rcl}
P(2n)&=&\left\{
\begin{array}{clcl}
2 &\hbox{if }\ |V^{(2k-1)}|< 2n \leq |U^{(2k-1)}|\ \hbox{ for some } k \in \mathbb N\,,\\
1 &\mbox{otherwise}.\\
\end{array}
\right.\\[4mm]
P(2n+1)&=&\left\{
\begin{array}{clcl}
3 & \hbox{if }\ |V^{(2k)}|< 2n+1 \leq |U^{(2k)}| \ \hbox{ for some } k \in \mathbb N\,,\\
2 &\mbox{otherwise}.\\
\end{array}
\right.
\end{array}
$$

\item[(ii)]
Let both $b$ and $a$ be even.
$$
\begin{array}{rcl}
P(2n)&=&\left\{
\begin{array}{clc}
2 & \hbox{if }\ |V^{(2k-1)}|< 2n \leq |U^{(2k)}| \hbox{ for some } k \in \mathbb N\,,\\
1 &\mbox{otherwise}.\\
\end{array}
\right.\\[4mm]
P(2n+1)&=&\left\{
\begin{array}{clc}
2 & \hbox{if }\quad 2n+1 \leq |U^{(1)}|=a-1\,,\\
2 & \hbox{if }\ |V^{(2k)}|< 2n+1 \leq |U^{(2k+1)}| \hbox{ for some } k \in \mathbb N\,,\\
1 &\mbox{otherwise}.\\
\end{array}
\right.
\end{array}
$$

\item[(iii)]
Let $b$ be odd and $a$ even.
$$
\begin{array}{rcl}
P(2n)&=&\left\{
\begin{array}{clc}
2 & \hbox{if }\ |V^{(3k-1)}|< 2n \leq |U^{(3k-1)}| \hbox{ for some } k \in \mathbb N\,,\\
1 &\mbox{otherwise}.\\
\end{array}
\right.\\[4mm]
P(2n+1)&=&\left\{
\begin{array}{clc}
3 & \hbox{if } |V^{(k)}|< 2n+1 \leq |U^{(k)}| \hbox{ for some } k \in \mathbb N,\
k\not\equiv 2\, {\rm mod}\, 3,\\
2 &\mbox{otherwise}.\\
\end{array}
\right.
\end{array}
$$

\item[(iv)]
Let both $b$ and $a$ be odd. We have
$$
\begin{array}{rcl}
P(2n)&=&\left\{
\begin{array}{clc}
1 & \hbox{if }\ 0\leq 2n \leq |U^{(1)}|=a-1,\\
0 &\mbox{otherwise}.\\
\end{array}
\right.\\[4mm]
P(2n+1)&=&\left\{
\begin{array}{clc}
2 & \hbox{if }\ 2n+1 \leq |V^{(1)}|=b,\\
4 & \hbox{if }\ |V^{(k)}| < 2n+1 \leq |U^{(k)}| \hbox{ for some } k \geq 2,\\
3 &\mbox{otherwise}.\\
\end{array}
\right.
\end{array}
$$
\end{enumerate}
\end{thm}

Note that we can either derive the values of the palindromic
complexity for both even and odd $n$ directly from
Proposition~\ref{allbidirect}, or we can determine only ${\cal
P}(2n)$ and then use the relation~\eqref{e:nerovnost} between
palindromic complexity and the first difference of factor
complexity,
$$
P(2n+1)=\triangle C(2n)+2-P(2n),
$$
knowing the first difference of factor complexity from
Theorem~\ref{t:faccom}.

\section{Conclusion}

The present paper completes the study of palindromic complexity of
infinite words associated with $\beta$-integers for $\beta$ with
eventually periodic R\' enyi expansion of unity in base $\beta$,
i.e., $\beta$ a Parry number. The study started already
in~\cite{AmFrMaPe1} for simple Parry numbers; here we focus on
non-simple Parry numbers $\beta$. We focus on quadratic non-simple
Parry numbers, since only in this case, one finds infinitely many
palindromes in $u_\beta$. Such infinite words were studied for
their balance properties in~\cite{balance}. Other properties to be
investigated are for example recurrent times or return words,
which are sofar described only for a class of simple Parry
numbers~\cite{return}.

\section*{Acknowledgements}

The authors acknowledge financial support by Czech Science
Foundation GA \v{C}R 201/05/0169, by the grant LC00602 of the
Ministry of Education, Youth, and Sports of the Czech Republic.


\end{document}